\documentclass[10pt]{article}
\usepackage{latexsym,amssymb,amsmath}

\def\N{\mathbb{N}}
\def\Q{\mathbb{Q}}
\def\Z{\mathbb{Z}}
\def\R{\mathbb{R}}
\def\C{\mathbb{C}}

\def\proof{\par\medskip\noindent{\em Proof. }}
\def\eproof{\hfill{$\Box$}\bigskip}
\def\sus{\subset}
\def\al{\alpha}
\def\be{\beta}

\def\de{\delta}

\def\ds{\dots}

\def\ord{\mathrm{ord\,}}

\title{Computing the $n$-th coefficient of an algebraic power series modulo $p$ in $O(\log n)$ operations}
\author{Martin Klazar\footnote{{\tt klazar@kam.mff.cuni.cz}}}

\begin{document}
\maketitle
\begin{abstract}
This is an exposition, for pedagogical purposes, of the formal power series proof of Bostan, Christol and Dumas \cite{bost_chri_duma}
of the result stated in the title (a corollary of the Christol theorem).
\end{abstract}

\begin{center}
{\bf Introduction}
\end{center}

The generating function $c=c(x)=\sum_{n\ge1}c_nx^n$ of the Catalan numbers 
$c_n=\frac{1}{n}\binom{2n-2}{n-1}=1,1,2,5,14,42,\ds$ satisfies the quadratic equation
$$
c=c^2+x\;.
$$
Can we effectively determine parity of $c_n$? Very easily: modulo $2$ one has
$$
c(x)=c(x)^2+x\equiv c(x^2)+x
$$
(by Lemma C below), and therefore $c_n\equiv c_{n/2}$ mod $2$ for $n>1$ and $c_1\equiv 1$ mod $2$ ($c_{n/2}=0$ 
if $n/2\not\in\N_0$). Hence 
\begin{quote}
$c_n\equiv 1$ mod $2$  for $n=2^m$, $m\in\N_0$, and $c_n\equiv 0$ mod $2$  else\,. 
\end{quote}
And $c_n$ modulo $3$ or modulo other prime $p$? An answer to this is provided much more generally by the 
{\em Christol theorem}, proved by Christol \cite{chri} for $q=2$ and by Christol, Kamae, Mend\` es France and Rauzy \cite{chri_al} 
for any prime power $q$ (see also Allouche and Shallit \cite[Theorem 12.2.5]{allo_shal}): a power series 
$f(x)\in F_q[[x]]$ is algebraic over $F_q(x)$ if and only if one can generate its coefficients $f_0,f_1,\ds$ 
by a DFAO (a deterministic finite automaton with output that reads the $l+1$ 
digits of the $q$-adic expansion $n=(n_l\ds n_1n_0)_q$ and outputs $f_n$, we give one example in Concluding remarks). It is immediate from the Christol theorem
that one can compute the $n$th coefficient of an algebraic power series in $F_q[[x]]$ in $O(\log n)$ arithmetic 
operations because $n$ has $O(\log n)$ $q$-adic digits.

This write-up is an exposition of the proof of this result, stated as Theorem A below, that was given by means 
of power series in the recent preprint of Bostan, Christol and Dumas \cite{bost_chri_duma}. I promised in my 
course `Kombinatorick\'e po\v c\'\i t\'an\'\i' (taught in this summer semester 2015/16) to supplement 
my unfortunately not very clear oral presentation of (something like) Theorem A with a written one --- here it is. 

\begin{center}
{\bf Algebraic power series in $F_p[[x]]$}
\end{center}

First we review some notation and notions from algebra. Puiseux series $K((x))_P$ and their order 
are only used at the end in Proposition H that is not needed to prove Theorem A.
$\N=\{1,2,\ds\}$, $\N_0=\{0,1,\ds\}$, $\Z$ is the ring of integers, $F_p=(\Z/p\Z,+,\cdot)$ 
is the finite $p$-element field where $p$ is a prime number, and $[x^n]f$ denotes the coefficient 
of $x^n$ in $f$. For $p\in\N$, $p\ge2$, the $p$-adic expansion $m=(m_l\ds m_1m_0)_p$ of $m\in\N_0$ 
is given by the unique expression 
$m=\sum_{i=0}^lm_ip^i$ where $m_i\in\{0,1,\ds,p-1\}$ and (for $m>0$) $m_l\ne0$; $0=(0)_p$. Let $K$ be a field or an
integral domain. For a field we denote by $\overline{K}$ its algebraic closure (e.g. $\overline{\R}=\C$ or 
$\overline{\Q}$ are the algebraic numbers). We denote by $K[[x]]\sus K((x))\sus K((x))_P$, respectively, the ring of 
power series and the fields of Laurent series and Puiseux series with coefficients in $K$. Their elements are 
the formal series (we treat here all 
power etc. series formally) $\sum_{n\ge k}a_{n/d}x^{n/d}$, $a_{n/d}\in K$, where $k=0$ and $d=1$ for p. s., $k\in\Z$ 
and $d=1$ for L. s., and $k\in\Z$ and $d\in\N$ for P. s. For $f=f(x)\in K((x))$ we use notation $f_n:=[x^n]f(x)$. 
From the algebraic point of view, $K((x))$ is the field of 
fractions of $K[[x]]$ and for $K$ of characteristic $0$ (e.g. $K=\Q$) $\overline{K}((x))_P$ is the algebraic closure
of $K((x))$. Thus in $\overline{K}((x))_P$ every polynomial with coefficients in $K((x))$ has a root. Recall that 
$K[x]\sus K(x)$ is, respectively, the ring of polynomials with coefficients in $K$ and its field of fractions ($K(x)$ 
are the rational `functions' over $K$). In fact, $K(x)$ is naturally a subfield of $K((x))$ and $K[x]$ a subring of 
$K[[x]]$. By $\deg a$, $a\in K[x]$, we denote the degree of a polynomial and by $\deg_y a$, $a\in K[x,y]$, the degree
$\deg a$ when $a$ is understood as $a\in K[x][y]$. By $\ord f$, $f\in K((x))_P$, we denote the order of $f$, the 
minimum $n/d$ such that $a_{n/d}\ne0$. We set $\deg 0=-\infty$ and $\ord 0=+\infty$. For any $f,g\in K((x))_P$ one has 
$\ord(fg)=\ord f+\ord g$ and $\ord(f+g)\ge\min(\ord f,\ord g)$, with equality if $\ord f\ne\ord g$. Similarly, for any 
$f,g\in K[x]$ one has $\deg(fg)=\deg f+\deg g$ and $\deg(f+g)\le\max(\deg f,\deg g)$, with equality if $\deg f\ne\deg g$.

We state the basic result on computing coefficients of an algebraic power series in $F_p[[x]]$ (this is essentially
proven in $\cite{chri_al}$ but the formulation below is ours). 

\medskip\noindent
{\bf Theorem A (\cite{chri_al}, 1980). }{\em Suppose that $p$ is prime, $P\in\Z[x,y]$ is a nonzero polynomial,
$f=f_0+f_1x+f_2x^2+\ds\in\Z[[x]]$ is a power series satisfying $P(x,f(x))=0$, and 
$$
h=(d+1)(p^d-d+1)\deg_xP\ \mbox{ where }\ d=\deg_yP\ (\ge1)\;. 
$$
Then from given $p,P,f_0,f_1,\ds,f_h$ one can construct an algorithm 
(actually a DFAO) that computes
$$
n\mapsto f_n\ \mathrm{mod}\ p
$$
in $O(\log n)$ arithmetic operations ($n\ge2$). 
}

\medskip\noindent
But which operations exactly? These are arithmetic operations ($+,-,\times,:$) in the field 
$F_p$ and the operation of division with remainder in the ring $\Z$. By `one can construct'  here and below we mean 
that the task can be performed constructively (i.e. the algorithms or objects not merely exist 
but exist effectively). Theorem A follows from the following more specific result (again the formulation is ours). 

\medskip\noindent
{\bf Theorem B (\cite{bost_chri_duma}). }{\em Let $p$, $P$, $f$, and $h$ be as in Theorem A. We denote the 
mod $p$ reduction of $f$ again by 
$f\in F_p[[x]]$. Then from given $p,P,f_0,f_1,\ds,f_h$ one can construct polynomials $a,c\in F_p[x]$, a number 
$e\in\N$, a row $u\in F_p^{1\times e}$, $p$ matrices $A_0,A_1,\ds,A_{p-1}\in F_p^{e\times e}$, and a column 
$v\in F_p^{e\times 1}$ such that for any $n\in\N_0$ one has
\begin{eqnarray*}
[x^n]f(x)&=&f_n=[x^n](c(x)+a(x)h(x))\ \mbox{ where }\\
&&h(x)=\sum_{m=0}^{\infty}(uA_{m_l}\ds A_{m_1}A_{m_0}v)x^m\in F_p[[x]]
\end{eqnarray*} 
($m=(m_l\ds m_1m_0)_p$ is the $p$-adic expansion of $m$).
}

\medskip
Polynomial $P$ in Theorem A in general does not determine unique power series solution $y=f(x)$ of 
$P(x,y)=0$, e.g. $y^2-x^2=0$ has two solutions $y=x$ and $y=-x$. One can introduce restrictions, for example to require 
that $P(0,0)=0$, $P_y(0,0)\ne0$, and $f(0)=0$, producing unique solution, but we want result for completely general $P$. So to determine $f(x)$ we need to know a few initial coefficients, which is 
captured by the parameter $h$ in Theorems A and B. Bounds on $h$ are considered in Corollary G and 
Proposition H. 

Let us see how Theorem B implies Theorem A. We are given a prime $p$, a nonzero $P\in\Z[x,y]$, 
and first few coefficients $f_0,f_1,\ds,f_h$ of a solution $y=f(x)\in\Z[[x]]$ of the equation $P(x,y)=0$, 
where $h$ is as in Theorem A. (We assume this is an honest input, $f_0,f_1,\ds,f_h$ are really initial 
coefficients of a solution. We leave as an exercise for the reader to design an effective check
that it is true.) Taking out  the largest common factor of the coefficients of $P(x,y)$ we may assume 
that they are together coprime. Then reducing $P(x,f(x))=0$ modulo $p$ we get the relation $E(x,f(x))=0$ 
where $E\in F_p[x,y]$ is the {\em nonzero} mod $p$ reduction of $P$ and $f\in F_p[[x]]$ 
denotes again the mod $p$ reduction of $f$ --- its coefficients $f_n$ we want to compute. 
We do it by Theorem B. We construct in $O_{p,P}(1)$ operations the objects $a,c,e,u,A_0,A_1,\ds,A_{p-1}$, and $v$. 
For given $n\in\N_0$ we compute $f_n$ as
$$
f_n=[x^n]c(x)+\sum_{m=n-\deg a}^n[x^{n-m}]a(x)\cdot uA_{m_l}\ds A_{m_1}A_{m_0}v
$$
where $m=(m_l\ds m_1m_0)_p$. For each $m$ ($>0$) in the summation range it takes 
$l+1=l(m)+1=\lfloor\log_pm\rfloor+1$ operations in $\Z$ (divisions by $p$ with remainder) to determine its $p$-adic 
expansion. It takes $l+1$ matrix times column multiplications, one row times column multiplication, and one 
scalar multiplication (in $F_p$) to evaluate the summand. So we evaluate $f_n$ in 
$$
\sum_{m=n-\deg a}^n(l(m)+1+(l(m)+1)e(2e-1)+2e+1)=O((\deg a)e^2\log_pn)
$$
operations. Altogether we need $$O_{p,P}(1)+O((\deg a)e^2\log_pn)=O_{p,P}(\log n)$$ operations 
($n\ge2$). 

\begin{center}
{\bf Proof of Theorem B}
\end{center}

We start with four lemmas. Lemma C is well known and is crucial for the existence of the algorithm of Theorem A, and 
in fact for the all algebra (and analysis) in characteristic $p$.

\medskip\noindent
{\bf Lemma C. }{\em If $p$ is prime and $z\in F_p((x))$ is any Laurent series then $z(x)^p=z(x^p)$. Consequently,
for any $k\in\N_0$ one has 
$$
z(x)^{p^k}=z(x^{p^k})\;.
$$
}
\vspace{-3mm}
\noindent
\proof
Let $z(x)=z_rx^r+z_{r+1}x^{r+1}+\cdots$. Since $a^p=a$ for any $a\in F_p$ and $(u+v)^p=u^p+v^p$ for any $u,v\in F_p((x))$
(because $\binom{p}{i}\equiv0$ mod $p$ for $0<i<p$), we indeed have 
\begin{eqnarray*}
z(x)^p&=&z_rx^{pr}+(z_{r+1}x^{r+1}+\ds)^p\\
&=&z_rx^{pr}+z_{r+1}x^{p(r+1)}+(z_{r+2}x^{r+2}+\ds)^p\\
&\vdots&\\
&=&z(x^p)\;.
\end{eqnarray*}
\eproof

We define operators $S_r:\;F_p[[x]]\to F_p[[x]]$, $r\in\{0,1,\ds,p-1\}$, by 
$${\textstyle
S_rz(x)=S_r(z(x))=S_r(\sum_{n\ge0}z_nx^n)=\sum_{n\ge0}z_{pn+r}x^n\;.}
$$ 
Thus $S_rz(x)$ arises from $z(x)$ by taking only terms $z_nx^n$ with $n\equiv r$ modulo $p$
and replacing the $n$ in the exponent with $(n-r)/p$.

\medskip\noindent
{\bf Lemma D (\cite{bost_chri_duma}). }{\em Operators $S_r$ have the following properties.
\begin{enumerate}
\item Linearity, $S_r(au+bv)=aS_r(u)+bS_r(v)$ for any $a,b\in F_p$ and $u,v\in F_p[[x]]$.
\item $S_r(uv)=\sum_{s+t\equiv r\;\mathrm{mod}\;p}x^{\lfloor(s+t)/p\rfloor}S_s(u)S_t(v)$ for any $u,v\in F_p[[x]]$.
\item $S_r(u(x)v(x^p))=S_r(u(x))v(x)$ for any $u,v\in F_p[[x]]$.
\item If $u\in F_p[[x]]$ and $n=(n_l\ds n_1n_0)_p\in\N_0$ then 
$$
u_n=[x^n]u(x)=(S_{n_l}\ds S_{n_1}S_{n_0}u(x))(0)=[x^0]S_{n_l}\ds S_{n_1}S_{n_0}u(x)\;.
$$
\item For any $0\ne u\in F_p[[x]]$ there is an $r\in\{0,1,\ds,p-1\}$ such that $S_ru\ne0$.
\item $S_rF_p[x]\sus F_p[x]$ and $\deg S_ra\le(\deg a)/p$ for any $a\in F_p[x]$.
\end{enumerate}
}

\noindent
\proof
1. This is immediate from the definition.

2. We have 
$$
[x^n]S_r(uv)=\sum_{i+j=pn+r}u_iv_j=\sum_{l+m=n,n-1;s+t=r,r+p}u_{pl+s}v_{pm+t}
$$ 
($i,j,l,m,n\in\N_0$, $r,s,t\in\{0,1,\ds,p-1\}$). The last sum equals the coefficient of $x^n$ in the sum stated in 2.

3. This is a corollary of 2: $S_t(v(x^p))$ is the zero power series for $t\ne0$ and $S_0(v(x^p))=v(x)$. 

4. Since ($r,s,r_i\in\{0,1,\ds,p-1\}$)
$$
S_ru(x)=\sum_{n\ge0}u_{pn+r}x^n,\ S_sS_ru(x)=\sum_{n\ge0}u_{p(pn+s)+r}x^n=\sum_{n\ge0}u_{p^2n+ps+r}x^n
$$
and so on, the constant term of $S_{r_l}\ds S_{r_1}S_{r_0}u(x)$ is just $u_{(r_l\ds r_1r_0)_p}$.

5. By the definition of $S_r$, for $r$ one may take the least significant $p$-adic digit of any $n\in\N_0$ for 
which $u_n=[x^n]u\ne0$.

6. Clear from the definition of $S_r$. 
\eproof

\medskip\noindent
{\bf Lemma E. }{\em Let $K$ be any field. If $v_1,\ds,v_{d+1}\in K[x]^d$ are $d+1$ $d$-tuples whose entries 
have degrees at most $c\in\N_0$, then one can construct coefficients $d_1,\ds,d_{d+1}\in K[x]$, not all $0$, such that
$$
\sum_{i=1}^{d+1}d_iv_i=(0,0,\ds,0)
$$
and $\deg d_i\le dc$ for every $i=1,\ds,d+1$.
}

\noindent
\proof
Let $M\in K[x]^{(d+1)\times d}$ be the matrix with rows $v_1,\ds,v_{d+1}$. We may assume that it has
full rank $d$, its columns are linearly independent over $K(x)$. (If not, we find by standard linear algebra some
maximal index set $I\sus\{1,2,\ds,d\}$ of linearly independent columns and solve the problem for restrictions 
of the tuples $v_i$ to coordinates with indices in $I$). We find $d$ row indices $I\sus\{1,2,\ds,d+1\}$ such that 
$\det M_I\ne0$, where $M_I$ arises from $M$ by deleting the row $j$ not in $I$ ($[d+1]\backslash I=\{j\}$). The displayed equation can be written as (we replace $d_i$ by $x_i$)
$$
(x_1,\ds,x_{j-1},x_{j+1},\ds,x_{d+1})M_I=-x_jv_j\;.
$$ 
We set $x_j=1$ and solve this linear system for the $x_i$, $i\in I$. 
By Cramer's rule, the solution is $x_i=\pm\det N_i/\det M_I$, $i\in I$, where $N_i$ arises from $M_I$ by replacing 
the row $i$ with the row $-x_jv_j=-v_j$ (the determinants are in $K[x]$). From the definition of determinant and bound on 
degrees of the entries in $M_I$ and $N_i$ we get that $\deg\det N_i,\deg\det M_I\le dc$. Thus the desired 
coefficients are $d_j=\det M_I$ and $d_i=\pm\det N_i$, $i\in I$. Then $\deg d_i\le dc$ for every $i=1,\ds,d+1$ and 
$d_j\ne0$.
\eproof

\medskip\noindent
{\bf Lemma F. }{\em Let $p$ be prime, $0\ne E\in F_p[x,y]$, $d=\deg_yE$, and let $f\in F_p[[x]]$ satisfy $E(x,f(x))=0$ (so $d\ge1$). Then from given polynomial $E$ one can construct polynomials $c_0,c_1,\ds,c_k\in F_p[x]$,
$0\le k\le d$, such that we have equation
$$
\sum_{i=0}^kc_i(x)f(x^{p^i})=0\;,
$$
$c_0\ne0$, and $\deg c_0\le(d+1)(p^d-d+1)\deg_xE$.
}

\noindent
\proof
We have $a_df^d=\sum_{i=0}^{d-1}a_if^i$ where $a_i=[y^i]E\in F_p[x]$ and $a_d\ne0$. It follows that for given $E$, $0\le i\le d$, and $n\in\N_0$ one can construct polynomials $a_{i,n}\in F_p[x]$ such that 
$${\textstyle
a_{d,n}f^n=\sum_{i=0}^{d-1}a_{i,n}f^i,\ a_{d,n}\ne0\;.
}
$$ 
Thus the $d+1$ elements $f^{p^i}=f(x^{p^i})$ (Lemma C), $i=0,1,\ds,d$, are linearly dependent over $F_p(x)$ and one 
can construct the dependency relation
$$
{\textstyle
\sum_{i=0}^dd_i(x)f(x^{p^i})=0
}
$$
with $d_i\in F_p[x]$ that are not all zero. Let $j\in\N_0$ be minimum
with $d_j\ne0$. Using Lemma D.5, we can construct $r_1,\ds,r_{j-1}\in\{0,1,\ds,p-1\}$ such that  
$S_{r_1}\ds S_{r_{j-1}}d_j(x)\ne0$. Applying $S_{r_1}\ds S_{r_{j-1}}$ on the dependency relation and using Lemma D.3 
we get the desired equation
$$
{\textstyle
\sum_{i=j}^d(S_{r_1}\ds S_{r_{j-1}}d_i(x))f(x^{p^{i-j}})=0\;.
}
$$
We set $c_i=S_{r_1}\ds S_{r_{j-1}}d_{i+j}(x)$, $i=0,1,\ds,k=d-j$, and have $c_0=S_{r_1}\ds S_{r_{j-1}}d_j(x)\ne0$.

We derive the upper bound on $\deg c_0$. Since $S_r$ does not increase degree, $\deg c_0\le\max_i\deg d_i$. 
We have this recurrence for the polynomials $a_{i,n}$:
$a_{d,n}=1$ and $a_{i,n}=\de_{i,n}$ for $0\le i,n\le d-1$, $a_{i,d}=a_i$ for $0\le i\le d$, and, 
for $n>d$, $a_{d,n}=a_{d,n-1}a_d=a_d^{n-d+1}$ and $a_{i,n}=a_{d-1,n-1}a_i+a_{i-1,n-1}a_d$ for $0\le i\le d-1$, with $a_{-1,n-1}=0$.
Clearly, $\max_{0\le i\le d}\deg a_i=\deg_xE$. Using the recurrence and induction we get the bound 
$\max_i\deg a_{i,n}=0$ for $n=0,1,\ds,d-1$, $\max_i\deg a_{i,d}=\deg_xE$, and 
$\max_i\deg a_{i,n}\le(n-d+1)\deg_xE$ for $n\ge d$. Applying Lemma E on the $d+1$ $d$-tuples 
$(a_{i,n}\;|\;0\le i\le d-1)$, $n=1,p,p^2,\ds,p^d$, we get $d_i\in F_p[x]$ in the dependency relation with 
$\max_i\deg d_i\le(d+1)(p^d-d+1)\deg_xE$.
\eproof

\noindent
The lemma is inspired by Lemma 12.2.3 in Allouche and Shallit \cite{allo_shal}.
Equations of this type are called Mahler equations, after the German--British mathematician Kurt Mahler 
(1903--1988) who used them to prove transcendence of various numbers. 

We start the proof of Theorem B proper. Let $0\ne E\in F_p[x,y]$ be the mod $p$ reduction of $P(x,y)$ (with coprime
coefficients) and $d=\deg_yE\ge1$. Since $E(x,f)=0$, by Lemma F one can construct a Mahler equation 
$$
c_0(x)f(x)+c_1(x)f(x^p)+\ds+c_k(x)f(x^{p^k})=0
$$ 
for $f(x)$ where $c_i\in F_p[x]$, $k\le d$, $c_0(x)\ne0$, and $\deg c_0\le(d+1)(p^d-d+1)\deg_xE\le h$ ($h$ is 
as in Theorem A). We change the variable $f(x)$ to $g(x)$ by $f(x)=c_0(x)g(x)$. Using Lemma C and dividing by 
$c_0(x)^2$, we get a Mahler equation for $g(x)$, 
$$
g(x)+c_1(x)c_0(x)^{p-2}g(x^p)+\ds+c_k(x)c_0(x)^{p^k-2}g(x^{p^k})=0
$$ 
--- now the new coefficient $c_0$ is $1$. The small price we paid is that in general $g(x)$ is not a power series 
but $g(x)\in F_p((x))$. We write
$$
g(x)=\sum_{n<0}g_nx^n+\sum_{n\ge0}g_nx^n=:g_-(x)+h(x)\;,
$$
$g_-\in x^{-1}F_p[x^{-1}]$ and $h\in F_p[[x]]$ (this is the $h(x)$ of Theorem B). Denoting the Mahler equation for 
$g(x)$ as $L(x,M)g(x)=0$ (here $L(x,M)$ is a skew polynomial in $x$ and the operator $M:\;u(x)\mapsto u(x)^p=u(x^p)$) 
we get effectively a Mahler equation for the power series $h(x)$,
$$
L(x,M)h(x)=-L(x,M)g_-(x)=:b(x)
$$
or more explicitly
$$
h(x)+c_1(x)c_0(x)^{p-2}h(x^p)+\ds+c_k(x)c_0(x)^{p^k-2}h(x^{p^k})=b(x)\;.
$$
Necessarily $b(x)\in F_p[x]$ because $b(x)=-L(x,M)g_-(x)\in F_p[x^{-1},x]$ but also 
$b(x)=L(x,M)h(x)\in F_p[[x]]$. We write the equation as 
$$
h(x)=b(x)+d_1(x)h(x^p)+\ds+d_k(x)h(x^{p^k})
$$
where $b(x),d_i(x):=-c_i(x)c_0(x)^{p^i-2}\in F_p[x]$. Let $D:=\max(\deg b,\deg d_i,1\le i\le k)$ and $h_{(i)}:=h(x^{p^i})$. We consider the vector space 
$$
V:=\{\al+\be_0h_{(0)}+\ds+\be_kh_{(k)}\;|\;\al,\be_i\in F_p[x],\deg\al,\deg\be_i\le D\}
$$
over $F_p$. Clearly, $h(x)=h_{(0)}\in V$. It follows that $S_r(V)\sus V$ for each of the $p$ operators $S_r$: 
using the equation for $h(x)$ and Lemma D.3 (and D.1) we have
\begin{eqnarray*}
&&S_r(\al+\be_0h_{(0)}+\ds+\be_kh_{(k)})={\textstyle S_r(\al+\be_0b+\sum_{i=1}^k(\be_0d_i+\be_i)h_{(i)})
}\\
&&={\textstyle S_r(\al+\be_0b)+\sum_{i=1}^kS_r(\be_0d_i+\be_i)h_{(i-1)}
}
\end{eqnarray*}
and $\deg S_r(\al+\be_0b),\deg S_r(\be_0d_i+\be_i)\le \deg(\al+\be_0b)/p,\deg(\be_0d_i+\be_i)/p\le 2D/p\le D$ 
(Lemma D.6). The $(D+1)(k+2)$-element set
$$
B:=\{x^j,x^jh_{(i)}=x^jh(x^{p^i})\;|\;0\le j\le D,0\le i\le k\}
$$
is a generating set for $V$ --- the linear span of $B$ over $F_p$ is $V$ (but $B$ is not a basis for $V$, as 
stated in \cite{bost_chri_duma}; $h_{(0)}$ is a linear combination of other elements of $B$).
We set $e:=(D+1)(k+2)$ and by fixing an ordering of $B$ make it a tuple $B=(b_1,\ds,b_e)$. We set 
$v\in F_p^{e\times 1}$ to be the column of coordinates of $h(x)$ in $B$ (it has all $0$s and one $1$), 
$A_r\in F_p^{e\times e}$ for $r=0,1,\ds,p-1$ to be the matrices of the operators $S_r$ with respect to $B$ (the $j$-th column of $A_r$ lists the coordinates of $S_r(b_j)$ in $B$), and $u:=(b_1(0),\ds,b_e(0))\in F_p^{1\times e}$ to be the row of constant terms in $B$. Clearly we can obtain $v,A_r$, and $u$ effectively. By Lemma D.4 and the definition of $v,A_r$, 
and $u$ we have for any $n=(n_l\ds n_1n_0)_p\in\N_0$ that 
$$
h_n=[x^n]h(x)=[x^0]S_{n_l}\ds S_{n_1}S_{n_0}h(x)=uA_{n_l}\ds A_{n_1}A_{n_0}v\;.
$$
Returning to $f(x)$ we see that
$$
f(x)=c_0(x)g(x)=c_0(x)g_-(x)+c_0(x)h(x)
$$
and set, finally, $c(x):=c_0(x)g_-(x)$ and $a(x):=c_0(x)$. It follows that $a,c\in F_p[x]$. How do we get $g_-(x)$? 
It is the negative part of $g(x)=f(x)/c_0(x)$. Since $\deg c_0\le h$, we compute $g_-(x)$ from $c_0(x)$ and 
$f_0,f_1,\ds,f_h$.

Let us review the computation. We are given $p,P$, and $f_0,f_1,\ds,f_h$ (initial coefficients of a solution 
$y=f(x)\in\Z[[x]]$ to $P(x,y)=0$, $h$ is as given in Theorem A). Reducing $P$ modulo $p$ we get a nonzero $E\in F_p[x,y]$. From $E$ we compute 
by Lemma F the coefficients $c_0,\ds,c_k\in F_p[x]$, $c_0\ne0$, $\deg c_0\le h$, and $k\le\deg_yE$, of a Mahler 
equation for $f(x)$. We compute $g_-(x)\in x^{-1}F_p[x^{-1}]$ as the negative part of $(f_0+f_1x+\ds+f_hx^h)/c_0(x)$. 
Thus we have $L(x,M)=M^0+\sum_{i=1}^kc_i(x)c_0(x)^{p^i-2}M^i$ (the Mahler equation for $g(x)$) and compute 
$b(x)=-L(x,M)g_-(x)\in F_p[x]$. We have $K(x,M):=-\sum_{i=1}^kc_i(x)c_0(x)^{p^i-2}M^i$ (the relation 
$h(x)=b(x)+K(x,M)h(x)$) and set $D=\max(\deg b,\deg_xK(x,M))$ and $e=(D+1)(k+2)$. We take the $e$-tuple 
$$
B=(x^j,x^jh(x^{p^i})\;|\;0\le j\le D,0\le i\le k)
$$ and compute the column $v$ of coordinates of $h(x)$ in $B$ 
(this is easy), the matrices $A_0,\ds,A_{p-1}$ of the operators $S_r$ with respect to $B$ ($S_r$ acts on $x^j$ in the clear way, 
$S_r(x^jh(x^{p^i}))=S_r(x^j)h(x^{p^{i-1}})$ for $i>0$, and for $i=0$ we replace $h(x)$ with $b(x)+K(x,M)h(x)$), 
and the row $u$ of constant terms in $B$ (note that the constant term $h(0)$ is the constant term in 
$(f_0+f_1x+\ds+f_hx^h)/c_0(x)$). Finally, we compute $c(x)=c_0(x)g_-(x),a(x)=c_0(x)\in F_p[x]$. Thus we have computed 
$a,c,e,u,A_0,\ds,A_{p-1},v$, all we need to quickly evaluate $f_n$ modulo $p$. We are done. 
\eproof
\vspace{-3mm}
\begin{center}
{\bf Concluding remarks}
\end{center}

An interesting problem is to locate the first different coefficient of two different power series solutions of the same
polynomial equation. Theorem B yields the following corollary.

\medskip\noindent
{\bf Corollary G. }{\em If $p$ is prime, $0\ne E\in F_p[x,y]$, and $f,g\in F_p[[x]]$, $f\ne g$, satisfy 
$E(x,f(x))=E(x,g(x))=0$, then
$$
f_n\ne g_n\ \mbox{ for some }\ n\le (d+1)(p^d-d+1)\deg_xE,\ d=\deg_yE\;. 
$$
}

\vspace{-3mm}
\noindent
Another bound is given in  

\medskip\noindent
{\bf Proposition H. }{\em If $H$ is any field of characteristic $0$, $0\ne P\in H[x,y]$, and $f,g\in H[[x]]$, $f\ne g$, satisfy $P(x,f(x))=P(x,g(x))=0$, then
$$
f_n\ne g_n\ \mbox{ for some }\ n\le\frac{d^2+d-4}{2}\deg_xP,\ d=\deg_yP\;.
$$}
{\vskip-5mm}
\proof
We work with the fields $K=H(x)\sus H((x))=L$ and the algebraic closure $M=\overline{L}=\overline{H}((x))_P$, 
the field of Puiseux series with coefficients in $\overline{H}$. There is a nonzero polynomial 
$Q\in H[x,y]$ such that $Q$ divides $P$ in $H[x,y]$, $Q(x,f(x))=Q(x,g(x))=0$, 
and $Q$ as $Q(y)\in K[y]$ has in $M$ only simple roots. Indeed, write $P(y)=P_1(y)^{a_1}P_2(y)^{a_2}\ds P_k(y)^{a_k}$ 
where each $P_i\in K[y]$ is irreducible (in $K[y]$), $a_i\in\N$, and no $P_i$ divides other $P_j$, $i\ne j$. Irreducibility of 
the $P_i(y)$s implies as usual that they have in $M$ only simple roots and do not share roots. By the Gauss lemma 
(Lang \cite[Chapter IV.2]{lang}) we may take $P_i\in H[x][y]$. Hence we have $Q=P_i(y)$ where $P_i(f)=P_i(g)=0$ or $Q=P_iP_j$, $i\ne j$, where 
$P_i(f)=P_j(g)=0$.

Let $f_1=f,f_2=g,\ds,f_d\in M$ be the (different) roots of $Q(y)=a_dy^d+\ds+a_1y+a_0$, $a_i\in H[x]$ and 
$a_d\ne0$, so $d=\deg_yQ\in\N$ (and $d\ge2$). Let $e=\max_i\deg a_i=\deg_xQ\in\N_0$. From Vi\`eta's formula 
$(-1)^ia_{d-i}/a_d=\sum_{1\le j_1<\ds<j_i\le d}f_{j_1}\ds f_{j_i}$ and properties of the $\ord$ function we get
$$
\min_{1\le i\le d}\ord f_i\ge\min_{0\le i\le d}\ord(a_i/a_d)\ge-e\;.
$$
We look at the order of the discriminant $D\in H[x]$ of $Q(y)$ (Lang \cite[Chapter IV.6]{lang}),
$$
D=a^{2d-2}_d\prod_{1\le i<j\le d}(f_j-f_i)^2=:(f-g)^2E\;.
$$
Since $Q(y)$ has no multiple root, $D$ is a nonzero homogeneous polynomial in $a_0,a_1,\ds,a_d$ with degree $2d-2$
and coefficients in $\Z$. So
\begin{eqnarray*}
2\,\ord(f-g)&=&\ord D-\ord\,E\le\deg D-\ord E\\
&\le& e(2d-2)-(2d-2)0-2(-e)(d(d-1)/2-1)\\
&=&e(d^2+d-4)
\end{eqnarray*}
and, since $d\le\deg_yP$ and $e\le\deg_xP$, the stated bound follows.
\eproof

\noindent
Alternatively, to determine uniquely an algebraic power series $f(x)$, instead of taking first few coefficients 
$f_0,f_1,\ds,f_h$ and the polynomial equation $f(x)$ satisfies, one can represent $f(x)$ 
by a system of polynomial equations such that $f(x)$ is the first coordinate in the tuple of unique solutions to 
the system. Then one can effectively do some operations with algebraic power series in such representation (which 
applies in fact to multivariate power series), see Alonso, Castro-Jim\'enez and Hauser \cite{alon_cast_haus} 
for this interesting topic. 

So how do the Catalan numbers $c_n$ behave modulo $3$?  According to this DFAO (taken from 
the note of Rowland \cite{rowl}):
\begin{eqnarray*}
&&a_101b_1,\,a_12e_2,\,b_10b_1,\,b_11c_2,\,b_12d_0,\,c_20c_2,\,c_21b_1,\,c_22d_0,\,d_0012d_0,\,e_20c_2,\\
&&e_21d_0,\,e_22e_2\;.
\end{eqnarray*}
It has five states $a_1,b_1,c_2,d_0$, and $e_2$, with the output mod $3$ residue in the index, and 
$5\cdot3=15$ transitions, with the input ternary digits written between the states. Computation starts 
always at $a_1$ and follows the ternary digits $n-1=(t_l\ds t_1t_0)_3$, read from the least significant $t_0$. 
For example, $6-1=(12)_3$ sends us in two steps from $a_1$ to $d_0$, and indeed $c_6=42$ is $0$ mod $3$. 

For a bound on the number of states of the DFAO obtained from an algebraic power series see Bridy \cite{brid}, 
and for further $f\in\Z[[x]]$ whose reduction $f_n$ mod $p$ (or mod $p^k$) can be computed by a DFAO see Rowland 
and Yassawi \cite{rowl_yass} (and the references therein). 

Finally, we quote the two main results of Bostan, Christol and Dumas \cite{bost_chri_duma} which quantify the dependence 
of the complexity bound on $P$ and $p$. Theorem 5 in \cite{bost_chri_duma} states: 
\begin{quote} 
Let $E$ be a polynomial in $F_p[x,y]_{h,d}$ such that $E(0,0)=0$ and $E_y(0,0)\ne0$, and let $f\in F_p[[x]]$ be its
unique root with $f(0)=0$. Algorithm 1 computes the $N$th coefficient $f_N$ of $f$ using 
$\tilde{O}(d^3h^2p^{3d}\log N)$ operations in $F_p$.
\end{quote}
Here $h$ and $d$ bound the $x$- and $y$-degree, respectively, and $\tilde{O}(\cdot)$ is $O(\cdot)$ with polylogarithmic
factors omitted. Proof of this theorem we (roughly) followed in our expose. Another algorithm, considerably more efficient, is presented in 
Theorem 11 in \cite{bost_chri_duma}:
\begin{quote} 
Let $E$ in $F_p[x,y]_{h,d}$ satisfy $E(0,0)=0$ and $E_y(0,0)\ne0$, and let $f\in F_p[[x]]$ be its
unique root with $f(0)=0$. One can compute the coefficient $f_N$ of $f$ in
$h^2(d+h)^2\log N+\tilde{O}(h(d+h)^5p)$ operations in $F_p$.
\end{quote}

\medskip\noindent
{\sc Charles University, KAM MFF UK, Malostransk\'e n\'am. 25, 118 00 Praha, Czechia}

\end{document}